\documentclass[12pt]{amsart}
\usepackage[OT2,T1]{fontenc}
\usepackage[russian,USenglish]{babel}

\usepackage{amssymb}
\usepackage{bbm}
\usepackage{mathrsfs}
\usepackage{xypic}

\makeatletter
\def\hsmash{\relax 
  \ifmmode\def\next{\mathpalette\mathhsm@sh}\else\let\next\makehsm@sh
  \fi\next}
\def\makehsm@sh#1{\setbox\z@\hbox{#1}\finhsm@sh}
\def\mathhsm@sh#1#2{\setbox\z@\hbox{$\m@th#1{#2}$}\finhsm@sh}
\def\finhsm@sh{\wd\z@\z@ \box\z@}
\makeatother

\newtheorem{fac}{Fact}[section]
\newtheorem{subl}[fac]{Sublemma}
\newtheorem{lem}[fac]{Lemma}
\newtheorem{prop}[fac]{Proposition}
\newtheorem{theo}[fac]{Theorem}
\newtheorem{coro}[fac]{Corollary}

\theoremstyle{definition}

\theoremstyle{remark}
\newtheorem{rem}[fac]{Remark}

\newtheorem{ex}[fac]{Example}

\def\pmod#1{\nobreak\mkern8mu (\text{\rmfamily\upshape mod}\,\,#1)}
\def\cyr#1{\foreignlanguage{russian}{#1}}

\newcommand{\br}{ }
\newcommand{\brr}{, }

\newcommand{\Spec}{\mathop{\text{\rm Spec}}\nolimits}
\newcommand{\Div}{\mathop{\text{\rm Div}}\nolimits}
\newcommand{\Pic}{\mathop{\text{\rm Pic}}\nolimits}
\newcommand{\Gal}{\mathop{\text{\rm Gal}}\nolimits}
\newcommand{\PGL}{\mathop{\text{\rm PGL}}\nolimits}
\newcommand{\Aut}{\mathop{\text{\rm Aut}}\nolimits}
\newcommand{\Gr}{\mathop{\text{\rm Gr}}\nolimits}
\newcommand{\Bl}{\mathop{\text{\rm Bl}}\nolimits}
\newcommand{\Br}{\mathop{\text{\rm Br}}\nolimits}

\newcommand{\pr}{\mathop{\text{\rm pr}}\nolimits}
\newcommand{\ev}{\mathop{\text{\rm ev}}\nolimits}

\renewcommand{\div}{\mathop{\text{\rm div}}\nolimits}

\newcommand{\bbA}{{\mathbbm A}}
\newcommand{\bbC}{{\mathbbm C}}
\newcommand{\bbF}{{\mathbbm F}}
\newcommand{\bbN}{{\mathbbm N}}
\newcommand{\bbQ}{{\mathbbm Q}}
\newcommand{\bbR}{{\mathbbm R}}
\newcommand{\bbZ}{{\mathbbm Z}}

\newcommand{\calA}{{\mathscr{A}}}
\newcommand{\calK}{{\mathscr{K}}}
\newcommand{\calM}{{\mathscr{M}}}
\newcommand{\calO}{{\mathscr{O}}}
\newcommand{\calS}{{\mathscr{S}}}
\newcommand{\calU}{{\mathscr{U}}}

\newcommand{\calHC}{\mathscr{HC}}

\newcommand{\frakl}{\mathfrak{l}}
\newcommand{\frakp}{\mathfrak{p}}
\newcommand{\frakq}{\mathfrak{q}}
\newcommand{\frakr}{\mathfrak{r}}
\newcommand{\frakL}{\mathfrak{L}}

\newcommand{\reg}{{{\text{\rm reg}}}}

\newcommand{\n}{{{\text{\rm n}}}}

\newcommand{\Ab}{{\text{\bf A}}}
\newcommand{\Pb}{{\text{\bf P}}}

\newcommand{\ratarrow}{$%
$\definemorphism{rat}\dashed\tip\notip%
\spreaddiagramcolumns{-12pt}%
\, \raisebox{0.12mm}{$-$} \!\!\diagram%
\rrat & 
\enddiagram\!\!$%
$}

\newcounter{abc}
\newenvironment{abc}{\begin{list}{\rm \alph{abc}) }%
{\usecounter{abc} \leftmargin=0.0pt \labelsep=0.0pt %
\listparindent=0.0pt \labelwidth=0.0pt \parsep=\smallskipamount %
\itemsep=0.0pt \topsep=0.0pt \partopsep=\smallskipamount}}{\end{list}}

\newcounter{iii}
\newenvironment{iii}{\begin{list}{\rm \roman{iii}) }%
{\usecounter{iii} \leftmargin=0.0pt \labelsep=0.0pt %
\listparindent=0.0pt \labelwidth=0.0pt \parsep=\smallskipamount%
 \itemsep=0.0pt \topsep=0.0pt \partopsep=\smallskipamount}}{\end{list}}

\def\rightend#1#2{{%
 \leavevmode\nobreak\hskip .5em plus 1fil
 \penalty600 \hskip 0pt plus -1filll
 \vadjust{}\nobreak\hskip 0pt plus 1filll%
 #1\parfillskip=#2\relax \par}}

\def\eop{\ifmmode\rule[-22pt]{0pt}{1pt}\ifinner\tag*{$\square$}\else\eqno{\square}\fi\else\rightend{$\square$}{0pt}\fi}

\title[Del Pezzo surfaces of degree four violating the Hasse principle]{Del Pezzo surfaces of degree four\\ violating the Hasse principle\\ are Zariski dense in the moduli scheme}

\begin{document}

\author{J\"org Jahnel}

\address{D\'epartement Mathematik\\ Universit\"at Siegen\\ Walter-Flex-Stra\ss e~3\\ D-57068~Sie\-gen\\ Germany}
\email{jahnel@mathematik.uni-siegen.de}
\urladdr{http://www.uni-math.gwdg.de/jahnel}

\author{Damaris Schindler}

\address{Hausdorff Center for Mathematics\\ Endenicher Allee 62\\ D-53115 Bonn\\
Germany}
\email{damaris.schindler@hausdorff-center.uni-bonn.de}
\urladdr{http://www.math.uni-bonn.de/people/dschindl}


\date{October~31,~2014.}

\keywords{Del Pezzo surface, Hasse principle, moduli scheme}

\subjclass[2010]{\mbox{Primary 11G35; Secondary 14G25, 14J26, 14J10}}

\begin{abstract}
We show that, over every number field, the degree four del Pezzo surfaces that violate the Hasse principle are Zariski dense in the moduli scheme.
\end{abstract}

\maketitle

\section{Introduction}

A del Pezzo surface is a smooth, proper algebraic
surface~$S$
over a field
$K$
with very ample anti-canonical sheaf
$\calK^{-1}$.
Over~an algebraically
closed field, every del Pezzo surface of degree
$d \leq 7$
is isomorphic to
$\Pb^2$,
blown up in
$(9-d)$
points in general position~\cite[The\-o\-rem~24.4.iii)]{Ma}.

According to the adjunction formula, a smooth complete intersection of two quadrics in
$\Pb^4$
is del Pezzo. The converse is true, as~well. For every del Pezzo surface of degree four, its anticanonical image is the complete intersection of two quadrics in
$\Pb^4$~\cite[Theorem~8.6.2]{Do}.\smallskip

Although del Pezzo surfaces over number fields are generally expected to have many rational points, they do not always fulfill weak~approximation. Even~the Hasse principle may~fail. The~first example of a degree four del Pezzo surface violating the Hasse principle has been devised by B.~Birch and Sir Peter Swinnerton-Dyer \mbox{\cite[Theorem~3]{BSD}}. It~is given
in~$\Pb^4_\bbQ$
by the equations
\begin{eqnarray*}
             T_0T_1 & = & T_2^2 - 5T_3^2 \, ,\\
(T_0+T_1)(T_0+2T_1) & = & T_2^2 - 5T_4^2 \, .
\end{eqnarray*}
Meanwhile,~more counterexamples to the Hasse principle have been constructed, see, e.g.,~\cite[Examples~15 and~16]{BBFL}. Only~quite recently, N.\,D.\,Q.~Nguyen \mbox{\cite[Theorem~1.1]{Ng}} proved that the degree four del Pezzo surface, given~by
\begin{eqnarray*}
                        T_0T_1 & = & T_2^2 - (64k^2+40k+5)T_3^2 \, ,\\
(T_0+(8k+1)T_1)(T_0+(8k+2)T_1) & = & T_2^2 - (64k^2+40k+5)T_4^2
\end{eqnarray*}
is a counterexample to the Hasse principle if
$k$
is an integer such that
$64k^2+40k+5$
is a prime~number. In~particular, under the assumption of Schinzel's hypothesis, this family contains infinitely many members that violate the Hasse principle.\medskip

In this paper, we prove that del Pezzo surfaces of degree four that are counterexamples to the Hasse principle are Zariski dense in the moduli~scheme. In~particular, we establish, for the first time unconditionally, that their number up to isomorphism is~infinite.
Although~certainly the case of the base
field~$\bbQ$
is of particular interest, we work over an arbitrary number
field~$K$.

Before we can state our main results, we need to recall some notation and facts about the coarse moduli scheme of degree four del Pezzo surfaces.\smallskip

For this we consider a del Pezzo surface
$X$
of degree four given as the zero set of two quinary quadrics
$$Q_1(T_0, \ldots, T_4) = Q_2(T_0, \ldots, T_4) = 0 \, .$$
The pencil $\smash{(uQ_1+vQ_2)_{(u:v) \in \Pb^1}}$
of quadrics defined by the forms
$Q_1$
and~$Q_2$
contains exactly five degenerate~elements. The~corresponding five values
$t_1, \ldots, t_5 \in \Pb^1(\overline{K})$
of~$t := (u:v)$
are uniquely determined by the surface
$X$,
up to permutation and the natural operation of
$\smash{\Aut(\Pb^1) \cong \PGL_2(\overline{K})}$.

Let~$\calU \subset (\Pb^1)^5$
be the Zariski open subset given by the condition that no two of the five components~coincide. Then~there is an~isomorphism
$$j\colon \calU\!/(S_5 \times \PGL_2) \stackrel{\cong}{\longrightarrow} \calM$$
to the coarse moduli
scheme~$\calM$
of degree four del Pezzo surfaces~\cite[Section~5]{HKT}.

The~quotient of
$\calU$
modulo
$S_5$
alone is the space of all binary quintics without multiple roots, up to multiplication by~constants. This~is part of the stable locus in the sense of Geometric Invariant Theory, which is formed by all quintics without roots of multiplicity
$\geq\! 3$~\cite[Proposition~4.1]{MFK}.\smallskip

Furthermore,~classical~invariant theory teaches that, for binary quintics, there are three fundamental invariants
$I_4$,
$I_8$,
and~$I_{12}$
of degrees
$4$,
$8$,
and
$12$,
respectively, that define an open embedding
$$\iota\colon \calU\!/(S_5 \times \PGL_2) \hookrightarrow \Pb(1,2,3)$$
into a weighted projective~plane. This~result is originally due to Ch.~Hermite~\cite[Section~VI]{He}, cf.\ \cite[Paragraphs 224--228]{Sa}. A~more recent treatment from a computational point of view is due to A.~Abdesselam~\cite{Ab}.\smallskip

Altogether, this yields an open embedding
$I\colon \calM \hookrightarrow \Pb(1,2,3)$.
More~generally, every family
$\pi\colon \calS \rightarrow B$
of degree four del Pezzo surfaces over a base scheme
$B$
induces a morphism
$$I_\pi=I \colon B \to \Pb(1,2,3) \, ,$$
which we call the {\em invariant map\/} associated
with~$\pi$.

\begin{rem}
There cannot be a fine moduli scheme for degree four del Pezzo surfaces, as every such
surface~$X$
has at least 16 automorphisms~\cite[Theorem~8.6.8]{Do}. (The statement of Theorem 8.6.8 in \cite{Do} contains a misprint, but it is clear from the proof that the described quotient group may be isomorphic to one of the listed groups or may be trivial).
\end{rem}

We can now state our first main result in the following form.

\begin{theo}\label{theo1}
Let\/~$K$
be any number field,
$U_\reg \subset \Gr(2,15)_K$
the open subset of the Gra\ss mann scheme that parametrizes degree four del Pezzo surfaces, and\/
$\calHC_K \subset U_\reg(K)$
be the set of all degree four del Pezzo surfaces
over\/~$K$
that are counterexamples to the Hasse~principle.\smallskip

\noindent
Then~the image of\/
$\calHC_K$
under the invariant map
$$\smash{I\colon U_\reg \longrightarrow \Pb(1,2,3)_K}$$
is Zariski~dense.
\end{theo}

In Theorem \ref{theo1} we identify the space
$S^2((K^5)^*)$
of all quinary quadratic forms with coefficients
in~$K$
with~$K^{15}$.
This~is clearly a non-canonical~isomorphism.
To~give an intersection of two quadrics
in~$\Pb^4_K$
is then equivalent to giving a
$K$-rational
plane through the origin
of~$K^{15}$,
i.e.\ a
$K$-rational
point on the Gra\ss mann scheme
$\Gr(2,15)_K$.
The~open subset
$U_\reg \subset \Gr(2,15)_K$
that parametrizes non-singular surfaces is isomorphic to the Hilbert scheme~\cite{Gr1} of del Pezzo surfaces of degree four
in~$\Pb^4_K$.
We~will not go into the details of this as they are not necessary for our~purposes.

\begin{rem}
An analogous result for cubic surfaces has recently been established by A.-S.~Elsenhans together with the first author~\cite{EJ}. Our~approach is partly inspired by the methods applied in the cubic surface~case. The concrete construction of del Pezzo surfaces of degree four that violate the Hasse principle is motivated by the work \cite{Ng} of N.\,D.\,Q.\ Nguyen.
\end{rem}

In~fact, more is true than stated in Theorem~\ref{theo1}. Our second main result seems to be a strengthening of the first one, but is, in fact, more or less~equivalent. Our~strategy will be to prove Theorem~\ref{theo1} first and then to deduce Theorem~\ref{theo2} from~it.

\begin{theo}\label{theo2}
Let\/~$K$
be any number field,
$U_\reg \subset \Gr(2,15)_K$
the open subset of the Gra\ss mann scheme that parametrizes degree four del Pezzo surfaces, and\/
$\calHC_K \subset U_\reg(K)$
be the set of all degree four del Pezzo surfaces
over\/~$K$
that are counterexamples to the Hasse~principle.\smallskip

\noindent
Then\/~$\calHC_K$
is Zariski dense in\/
$\Gr(2,15)_K$.
\end{theo}

\begin{rem}[Particular
$K3$~surfaces
that fail the Hasse principle]
In his article~\cite{Ng}, N.\,D.\,Q.\ Nguyen also provides families of
$K3$~surfaces
of degree eight that violate the Hasse~principle.
These~$K3$~surfaces
allow a morphism
$p\colon Y \rightarrow X$
that is generically 2:1 down to a degree four del Pezzo
surface~$X$
that fails the Hasse~principle.
Since~$X(K) = \emptyset$,
the existence of the morphism alone ensures that
$Y(K) = \emptyset$.

Nguyen's construction easily generalizes to our setting. One~has to intersect the cone
$CX \subset \Pb^5$
over the del Pezzo surface with a quadric that avoids the~cusp. The~intersection $Y$
is then a degree eight
$K3$~surface,
provided it is smooth, which it is generically according to Bertini's theorem.
Thus,~$Y$
is a counterexample to the Hasse principle provided it has an adelic~point.

For~$Y$,
the failure of the Hasse principle may be explained by the Brauer-Manin obstruction (cf.~Section~\ref{sec_Brauer} for details).
If~$\alpha \in \Br(X)$
explains the failure
for~$X$
then
$p^*\alpha$
does so
for~$Y$.

However,~the
$K3$~surfaces
obtained in this way do clearly not dominate the moduli space of degree eight
$K3$~surfaces.
Indeed,~the pull-back homomorphism
$p^*\colon \Pic(X_{\overline{K}}) \to \Pic(Y_{\overline{K}})$
doubles the intersection numbers and is, in particular,~injective. This~means that
$Y$
has geometric Picard rank at least six, while a general degree eight
$K3$~surface
is of geometric Picard rank~one.
\end{rem}

\section{A family of degree four del Pezzo surfaces}

We~consider the surface
$S := S^{(D;A,B)}$
over a
field~$K$,
given by the equations
\begin{eqnarray}
T_0T_1               & = & T_2^2 - DT_3^2 \, , \label{eq_eins} \\
(T_0+AT_1)(T_0+BT_1) & = & T_2^2 - DT_4^2      \label{eq_zwei}
\end{eqnarray}
for
$A,B,D \in K$.
We~will typically assume that
$D$
is not a square
in~$K$
and that
$S$
is non-singular. 
If~$S$
is non-singular, then
$S$
is a del Pezzo surface of degree~four.

\begin{prop}
\label{nonsing}
Let\/~$K$
be a field of characteristic\/
$\neq \!2$
and\/
$A,B,D \in K$.

\begin{abc}
\item
Then~the surface\/
$S^{(D;A,B)}$
is non-singular if and only if\/
$ABD\neq 0$,
$A\neq B$,
and\/
$A^2-2AB+B^2-2A-2B+1 \neq 0$.
\item
If\/~$D \neq 0$
then
\/
$S^{(D;A,B)}$
has not more than finitely many singular~points.
\end{abc}\smallskip

\noindent
{\bf Proof.}
{\em
a)
If
$D=0$
then
$S$
is a cone over a cone over four points
in~$\Pb^2$.
In~this case,
$S$
is singular, whether some of the four points coincide or~not.
Let~us suppose that
$D \neq 0$
from now~on.

A~point
$(t_0 : \ldots  :t_4) \in S(K)$
is singular if and only if the Jacobian matrix
$$
\left(
\begin{array}{ccccc}
       t_1     &        t_0      & -2t_2 & 2Dt_3 &  0    \\
 2t_0+(A+B)t_1 & (A+B)t_0+2ABt_1 & -2t_2 &  0    & 2Dt_4 
\end{array}
\right)
$$
is not of full~rank. In~particular, we immediately see that
$(0:1:0:0:0) \in S(K)$
is a singular point in the case that
$A=0$
or~$B=0$.
Thus,~we may assume that
$A \neq 0$
and
$B \neq 0$.

Furthermore,~if
$(t_0 : \ldots  :t_4) \in S(K)$
is singular then
$t_0^2 = ABt_1^2$
and
$t_2t_3 = t_2t_4 = t_3t_4 = 0$.
There~is clearly no point
on~$S$
that fulfills
$t_1 = 0$
together with these equations. Hence,~we may normalize the coordinates to
$t_1 = 1$,
i.e.~to
$t_0 = \pm\sqrt{AB}$,
and distinguish three cases.\smallskip

\noindent
{\em First case.}
$t_2=t_3=0$.

\noindent
Then
$t_0t_1=0$
by relation (\ref{eq_eins}), which is a~contradiction.\smallskip

\noindent
{\em Second case.}
$t_2=t_4=0$.

\noindent
Then
$t_0+At_1=0$ or $t_0+Bt_1=0$,
i.e.\
$\pm\sqrt{AB}+A=0$ or $\pm\sqrt{AB}+B=0$,
by the second equation. This immediately yields
$A=B$.\smallskip

\noindent
On~the other hand,
if~$A=B$
then
$\smash{((-A):1:0:\pm\sqrt{A/D}:0) \in S(\overline{K})}$
are singular~points.\smallskip

\noindent
{\em Third case.}
$t_3=t_4=0$.

\noindent
Then the relations (\ref{eq_eins}) and (\ref{eq_zwei}) together show that
$(t_0+At_1)(t_0+Bt_1) = t_0t_1$,
i.e.\
$$(\pm\sqrt{AB}+A)(\pm\sqrt{AB}+B) = \pm\sqrt{AB} \, .$$
This~equality clearly implies
$2AB \pm\sqrt{AB}(A+B-1) = 0$,
hence
$(A+B-1)^2 = 4AB$
and
$A^2-2AB+B^2-2A-2B+1 = 0$.

On the other hand, if
$A^2-2AB+B^2-2A-2B+1 = 0$
then
$\frac{1-A-B}2$
is a square root
of~$AB$
and a direct calculation shows that
$$\textstyle \big( \frac{1-A-B}2:1:\pm\sqrt{\frac{1-A-B}2}:0:0 \big)$$
are two singular points
on~$S^{(D;A,B)}$.\smallskip

\noindent
b)
Every singular point satisfies the relation
$t_0/t_1 = \pm\sqrt{AB}$.
Furthermore,~at least two of the coordinates
$t_2$,
$t_3$,
and~$t_4$
must~vanish. Together these~conditions define six lines in
$\Pb^4$,
which collapse to three in the case that
$AB = 0$.

If~there were infinitely many singular points then at least one of these lines would be entirely contained
in~$S$.
But~this is not the case, as, on each of the six lines, one equation of the~form
$$F(T_1) = T_2^2, \quad F(T_1) = -DT_3^2, \quad {\rm or} \quad F(T_1) = -DT_4^2$$
remains from the equations
of~$S$.
}
\eop
\end{prop}

\begin{rem}
Assume that
$D \in K$
is a non-square and that
$S^{(D;A,B)}$
is non-singular. Then~there is neither a
\mbox{$K$-rational}
point
$(t_0 \!:\! t_1 \!:\! t_2 \!:\! t_3 \!:\! t_4) \in S(K)$
such that
$t_0 = t_1 = 0$,
nor one such that
$t_0 + At_1 = t_0 + Bt_1 = 0$.
Indeed,~the first condition implies
$(t_0 \!:\! t_1 \!:\! t_2 \!:\! t_3 \!:\! t_4) = (0 \!:\! 0 \!:\! 0 \!:\! 0 \!:\! 1)$,
while, in view of
$A \neq B$,
the second one implies
$(t_0 \!:\! t_1 \!:\! t_2 \!:\! t_3 \!:\! t_4) = (0 \!:\! 0 \!:\! 0 \!:\! 1 \!:\! 0)$.
Both~points do not lie
on~$S$.
\end{rem}

\section{A class in the Grothendieck-Brauer group}
\label{sec_Brauer}

It~is a discovery of Yu.\,I.~Manin~\cite[\S47]{Ma} that a non-trivial element
$\alpha \in \Br(S)$
of the Gro\-then\-dieck-Brauer group \cite{Gr2}, \cite[Chapter~IV]{Mi} of a
variety~$S$
may cause a failure of the Hasse principle. Today,~this phenomenon is called the Brauer-Manin obstruction. Its~mechanism works as~follows.

Let
$K$
be a number field,
$\frakl \subset \calO_K$
a prime ideal, and
$K_\frakl$
be the corresponding completion. The~Grothendieck-Brauer group is a contravariant functor from the category of~schemes to the category of abelian groups. In~particular, for an arbitrary
scheme~$S$
and a
\mbox{$K_\frakl$-rational}
point
$x\colon \Spec K_\frakl \to S$,
there is a restriction~homomorphism
$x^*\colon \Br(S) \to \Br(\Spec K_\frakl) \cong \bbQ/\bbZ$.
For~a Brauer class
$\alpha \in \Br(X)$,
we call
$$\ev_{\alpha,\frakl}\colon S(K_\frakl) \longrightarrow \bbQ/\bbZ \, , \quad x \mapsto x^*(\alpha)$$
the local evaluation map, associated
to~$\alpha$.
Analogously, for
$\sigma\colon K \hookrightarrow \bbR$
a real prime, there is the local evaluation map
$\ev_{\alpha,\sigma}\colon S(K_\sigma) \to \frac12\bbZ/\bbZ$.

\begin{prop}[The Brauer-Manin obstruction to the Hasse principle]
\label{BM}
Let\/~$S$
be a projective variety over a number
field\/~$K$
and\/~$\alpha \in \Br(S)$
be a Brauer class.\smallskip

\noindent
For~every prime ideal\/
$\frakl \subset \calO_K$,
suppose that\/
$S(K_\frakl) \neq \emptyset$
and that the local evaluation map\/
$\ev_{\alpha,\frakl}$
is constant. Analogously, assume that, for every real prime\/
$\sigma\colon K \hookrightarrow \bbR$,
one has\/
$S(K_\sigma) \neq \emptyset$
and that the local evaluation map\/
$\ev_{\alpha,\sigma}$
is constant. 
Denote the values of\/
$\ev_{\alpha,\frakl}$
and\/
$\ev_{\alpha,\sigma}$
by\/
$e_\frakl$
and\/
$e_\sigma$,
respectively. 
If,~in this~situation,
$$\sum_{\frakl\subset\calO_K} \!\!e_\frakl + \!\!\!\sum_{\sigma\!\colon\!\! K \hookrightarrow \bbR} \!\!\!\!\!e_\sigma \neq \,0 \in \bbQ/\bbZ$$
then\/
$S$
is a counterexample to the Hasse~principle.\medskip

\noindent
{\bf Proof.}
{\em
The assumptions imply, in particular, that
$S$
is not the empty scheme. Consequently, there are
\mbox{$K_\tau$-rational}
points
on~$S$
for every complex prime\/
$\tau\colon K \hookrightarrow \bbC$.
The~Hasse principle would assert that
$S(K) \neq \emptyset$.

On~the other hand, by global class field theory~\cite[Section~10, Theorem~B]{Ta} one has a short exact sequence
$$0\rightarrow \Br(K)\rightarrow \bigoplus_{\nu} \Br(K_\nu) \rightarrow \bbQ/\bbZ\rightarrow 0 \, ,$$
where the direct sum is taken over all places
$\nu$
of the number field
$K$.
Assume that there is a point
$x\colon \Spec K \to S$.
Then
$x^*(\alpha) \in \Br(\Spec K)$
is a Brauer class that naturally maps to an element of
$\bigoplus_\frakl \Br(K_\frakl) \oplus \bigoplus_\sigma \!\Br(K_\sigma) \cong \bigoplus_\frakl \!\bbQ/\bbZ \oplus \bigoplus_\sigma \!\frac12\bbZ/\bbZ$
of a non-zero~sum, which is a contradiction to the exactness of the above sequence.
}
\eop
\end{prop}

\begin{prop}
\label{Brauerklasse}
Let\/~$K$
be any field of characteristic\/
$\neq \!2$
and\/
$A,B,D \in K \!\setminus\! \{0\}$
be arbitrary~elements. Set\/
$S := S^{(D;A,B)}$.
Suppose that\/
$D$
is a non-square and that\/
$S$
is non-singular. Put,~finally,
$\smash{L := K(\sqrt{D})}$.

\begin{abc}
\item
Then the quaternion algebra (see\/ \cite[Section 15.1]{Pi} for the notation)
$$\textstyle \calA := \big( L(S), \tau, \frac{T_0+AT_1}{T_0} \big)$$
over the function field\/
$K(S)$
extends to an Azumaya algebra over the whole
of\/~$S$.
Here,~by
$\tau \in \Gal(L(S)/K(S))$,
we denote the nontrivial~element.
\item
Assume that\/
$K$
is a number field and denote by\/
$\alpha \in \Br(S)$
the Brauer class, defined by the extension
of\/~$\calA$.
Let\/
$\frakl$
be any prime
of\/~$K$.
\begin{iii}
\item
Let\/~$(t_0 \!:\! t_1 \!:\! t_2 \!:\! t_3 \!:\! t_4) \in S(K_\frakl)$
be a point and assume that at least one of the quotients\/
$(t_0+At_1)/t_0$,
$(t_0+At_1)/t_1$,
$(t_0+Bt_1)/t_0$,
and\/
$(t_0+Bt_1)/t_1$
is properly defined and non-zero. Denote~that
by\/~$q$.
Then\/
$$\ev_{\alpha,\frakl}(t_0 \!:\! t_1 \!:\! t_2 \!:\! t_3 \!:\! t_4) =
\left\{
\begin{array}{cl}
0       & \text{~if\/~} (q,D)_\frakl = 1 \, , \\
\frac12 & \text{~if\/~} (q,D)_\frakl = -1 \, ,
\end{array}
\right.
$$
for\/
$(q,D)_\frakl$
the Hilbert~symbol.
\item
If\/~$\frakl$
is split in\/
$L$
then the local evaluation map\/
$\ev_{\alpha,\frakl}$
is constantly~zero.
\end{iii}
\end{abc}\smallskip

\noindent
{\bf Proof.}
{\em
a)
First of all,
$\calA$
is, by construction, a cyclic algebra of degree two. In~particular,
$\calA$
is simple~\cite[Section~15.1, Corollary~d]{Pi}.
Furthermore,~$\calA$
is obviously a central
\mbox{$K(S)$-algebra}.

To~prove the extendability assertion, it suffices to show that
$\calA$
extends as an Azumaya algebra over each valuation ring that corresponds to a prime divisor
on~$S$.
Indeed,~this is the classical Theorem of Auslander-Goldman for non-singular surfaces \cite[Proposition~7.4]{AG}, cf.~\cite[Chapter~IV, Theorem~2.16]{Mi}.

For~this, we observe that the principal divisor
$\div((T_0+AT_1)/T_0) \in \Div(S)$
is the norm of a divisor
on~$S_L$.
In~fact, it is the norm of the difference of two prime divisors, the conic, given by
$\smash{T_0+AT_1 = T_2 - \sqrt{D}T_4 = 0}$,
and the conic, given by
$\smash{T_0 = T_2 - \sqrt{D}T_3 = 0}$.
In~particular,
$\calA$
defines the zero element in
$H^2(\langle\sigma\rangle, \Div(S_L))$.
Under~such circumstances, the extendability
of~$\calA$
over the valuation ring corresponding to an arbitrary prime divisor
on~$S$
is worked out in
\cite[Paragraph 42.2]{Ma}.\smallskip

\noindent
b.i)
The~quotients
$$\textstyle \frac{T_0+AT_1}{T_0} \!/\! \frac{T_0+AT_1}{T_1} = \frac{T_2^2-DT_3^2}{T_0^2} \, , \;\;
\frac{T_0+BT_1}{T_0} \!/\! \frac{T_0+BT_1}{T_1} = \frac{T_2^2-DT_3^2}{T_0^2} \, ,
\;\, {\rm and} \;\;
\frac{T_0+AT_1}{T_0} \!/\! \frac{T_0+BT_1}{T_0} = \frac{T_2^2-DT_4^2}{(T_0+BT_1)^2}$$
are norms of rational~functions. Thus,~each of them defines the zero class in
$H^2(\langle\sigma\rangle, K(S_L)^*) \subseteq \Br K(S)$,
and hence
in~$\Br S$.
In~particular, the four expressions
$(T_0+AT_1)/T_0$,
$(T_0+AT_1)/T_1$,
$(T_0+BT_1)/T_0$,
and
$(T_0+BT_1)/T_1$
define the same Brauer~class.

The~general description of the evaluation map, given in \cite[Paragraph~45.2]{Ma} shows that
$\ev_{\alpha,\frakl}(t_0 \!:\! t_1 \!:\! t_2 \!:\! t_3 \!:\! t_4)$
is equal to
$0$
or
$\frac12$
depending on whether
$q$
is in the image of the norm map
$N_{L_\frakL/K_\frakl}\colon L_\frakL^* \to K_\frakl^*$,
or not, for
$\frakL$
a prime of
$L$
lying
above~$\frakl$.
This~is exactly what is tested by the Hilbert
symbol~$(q,D)_\frakl$.\smallskip

\noindent
ii)
If~$\frakl$
is split in
$L$
then the norm map
$\smash{N_{K(S_{L_\frakL})/K(S_{K_\frakl})}\colon K(S_{L_\frakL})^* \to K(S_{K_\frakl})^*}$
is surjective. In~particular,
$\smash{\frac{T_0+AT_1}{T_0}} \in K(S_{K_\frakl})^*$
is the norm of a rational function
on~$S_{L_\frakL}$.
Therefore, it defines the zero class in
$H^2(\langle\sigma\rangle, K(S_{L_\frakL})^*) \subseteq \Br K(S_{K_\frakl})$,
and thus
in~$\Br S_{K_\frakl}$.
Finally,~we observe that every
\mbox{$K_\frakl$-rational}
point
$x\colon \Spec K_\frakl \to S$
factors via
$S_{K_\frakl}$.
}
\eop
\end{prop}

Geometrically,~on a rank four quadric
in~$\Pb^4$,
there are two pencils of~planes. In~our situation, these are conjugate to each other under the operation
of~$\smash{\Gal(K(\sqrt{D})/K)}$.
The~equation
$T_0=0$
cuts two conjugate planes out of the quadric~(\ref{eq_eins}) and the same is true
for~$T_1 = 0$.
The~equations
$T_0+AT_1=0$
and~$T_0+BT_1=0$
each cut two conjugate planes out of~(\ref{eq_zwei}).

\begin{rem}
A.~V\'arilly-Alvarado and B.~Viray \cite[Theorem~5.3]{VAV} prove for a certain class of degree four del Pezzo surfaces that the Brauer-Manin obstruction is the only obstruction to the Hasse principle and to weak~approximation. Their~result is conditional under the assumption of Schinzel's hypothesis and the finiteness of Tate-Shafarevich groups of elliptic curves and based on ideas of O.~Wittenberg \cite[Th\'eor\`eme~1.1]{Wi}.
The~class considered in~\cite{VAV} includes our family~(\ref{eq_eins},\,\ref{eq_zwei}).
\end{rem}

One might formulate our strategy to prove
$S^{(D;A,B)}(K) = \emptyset$
for
$K$
a number field and particular choices of
$A$,
$B$,
and~$D$
in a more elementary way as~follows.

Suppose~that there is a point
$(t_0 \!:\! t_1 \!:\! t_2 \!:\! t_3 \!:\! t_4) \in S(K)$.
Then~$(t_0, t_1) \neq (0,0)$.
Among
$(t_0+At_1)/t_0$,
$(t_0+At_1)/t_1$,
$(t_0+Bt_1)/t_0$,
and
$(t_0+Bt_1)/t_1$,
consider an
expression~$q$
that is properly defined and non-zero.
Then~show that, for every prime
$\frakl$
of~$K$
including the Archimedean ones, but with the exception of exactly an odd number, the Hilbert symbol
$(q,D)_\frakl$
is equal
to~$1$.
Finally,~observe that such a behaviour contradicts the Hilbert reciprocity law~\cite[Chapter~VI, Theorem~8.1]{Ne}.

In~other words, the element
$q \in K_\frakl$
belongs to the image of the norm map
$N \colon L_\frakL \to K_\frakl$,
for~$\smash{L := K(\sqrt{D})}$
and
$\frakL$
a prime
of~$L$
lying
above~$\frakl$,
for all but an odd number of~primes. And~this is incompatible with \cite[Chapter~VI, Corollary~5.7]{Ne} or~\cite[Theorem~5.1 together with 6.3]{Ta}.

\section{Unramified primes}

\begin{lem}
\label{P2fuenfPunkte}
Let\/~$K$
be any field of characteristic\/
$\neq \!2$
and\/
$A,B,D \in K$
be elements such that\/
$D \neq 0$.
Then~the minimal resolution of singularities\/
$\widetilde{S}$
of\/~$S := S^{(D;A,B)}$
is geometrically isomorphic to\/
$\Pb^2$,
blown up in five~points.\medskip

\noindent
{\bf Proof.}
{\em
{\em First step.}
$S_{\overline{K}}$
is a rational~surface.

\noindent
For this, we observe that the quadric hypersurface defined by equation~(\ref{eq_eins}) has
$(0 \!:\! 0 \!:\! 0 \!:\! 0 \!:\! 1)$
as its only singular point, while the hypersurface defined by~(\ref{eq_zwei}) is regular at that~point. According~to \cite[Book~IV, Paragraph~XIII.11, Theorem~3]{HP}, this implies that
$S_{\overline{K}}$
is a rational~surface.\smallskip

\noindent
{\em Second step.}
$S_{\overline{K}}$
has only isolated singularities, each of which is a rational double~point.

\noindent
Proposition \ref{nonsing}.b) implies that
$S_{\overline{K}}$
has only isolated~singularities. Moreover,~we observe that
$S_{\overline{K}}$
contains only finitely many~lines. Indeed,~there are only finitely many lines through each singular point, since
$S_{\overline{K}}$
is not a~cone. On~the other hand, a line contained in
$\smash{S_{\overline{K},\reg}}$
has self-intersection number
$(-1)$
according to the adjunction formula, and is therefore~rigid.

Now~consider the projection
$S'$
of 
$\smash{S_{\overline{K}}}$
from a non-singular
$\smash{\overline{K}}$-rational
point~$p$
that does not lie on any of the~lines. By~construction,
$S'$
is a cubic surface
over~$\overline{K}$.
The~projection map
$\pr\colon S_{\overline{K}} \!\setminus\! \{p\} \to S'$
blows up the point
$p$
and is an isomorphism everywhere~else.
Indeed,
$\pr$
separates points since
$\pr(p_1) = \pr(p_2)$
implies that
$p$,
$p_1$,
and~$p_2$
are collinear, which enforces that the line through these three points must be contained in the quadrics defined by (\ref{eq_eins}) and (\ref{eq_zwei}), a
contradiction. The~same argument shows that
$\pr$
separates tangent~directions.

In~other words, the blowup
$\Bl_p(S_{\overline{K}})$
is a cubic~surface. Clearly,~it has as many singular points
as~$S_{\overline{K}}$.
In~particular,
$\Bl_p(S_{\overline{K}})$
is normal~\cite[Theorem~23.8]{Mt}. Moreover,~as 
$S_{\overline{K}}$
is rational,
$\Bl_p(S_{\overline{K}})$
is a rational~surface.

It~is well known that there are two kinds of normal cubic surfaces.
Either~$\Bl_p(S_{\overline{K}})$
is the cone over an elliptic curve or it belongs to one of the 21 cases having only double points being
$ADE$,
as listed in \cite[Table 9.1]{Do}. The~former case is impossible as this is not a rational~surface.
Further,~$ADE$-singularities
are rational~\cite[page~135]{Ar}.\smallskip

\noindent
{\em Third step.}
Conclusion.

\noindent
Let~now
$\smash{\pi\colon \widetilde{S} \to S}$
be the resolution~map. The~adjunction formula shows that
$\smash{K_{\widetilde{S}} = \pi^*i^*(-H) + E}$,
where $H$ is a hyperplane section and
$E$
a divisor
on~$\smash{\widetilde{S}}$
supported on the exceptional fibers
of~$\pi$.
But,~as the singularities of
$S$
are rational double points, one necessarily has
$E = 0$~\cite[Proposition 8.1.10]{Do}.

This~yields that
$C \!\cdot\! K_{\widetilde{S}} \leq 0$
for every curve
$C \subset \widetilde{S}$.
Moreover,~$\smash{K_{\widetilde{S}}^2 = [i^*(-H)]^2 = 4}$.
In~other words,
$\smash{\widetilde{S}}$
is a generalized del Pezzo surface~\cite{CT} of degree~five. By~an observation of Demazure \cite[Proposition~0.4]{CT}, this implies that
$\smash{\widetilde{S}}$
is geometrically isomorphic to
$\smash{\Pb^2}$,
blown up in five~points.
}
\eop
\end{lem}

\begin{coro}
\label{ratpoint}
Let\/~$\bbF_{\!\ell}$
be a finite field of characteristic\/
$\neq \!2$
and\/
$A,B,D \in \bbF_{\!\ell}$
such that\/
$D \neq 0$.
Then\/~$S := S^{(D;A,B)}$
has a regular\/
$\bbF_{\!\ell}$-rational
point.\medskip

\noindent
{\bf Proof.}
{\em
By~Lemma~\ref{P2fuenfPunkte}, the minimal resolution of singularities
$\smash{\widetilde{S}}$
of~$S$
is geometrically isomorphic to
$\Pb^2$,
blown up in five~points. In~such a situation, the Weil conjectures have been established by A.~Weil himself \cite[page~557]{We}, cf.~\cite[Theorem~27.1]{Ma}.

At~least one of the eigenvalues of Frobenius
on~$\smash{\Pic(\widetilde{S}_{\overline\bbF_{\!\ell}})}$
is equal
to~$(+1)$.
Say,~the number of eigenvalues
$(+1)$
is exactly
$n \geq 1$.
The~remaining
$(6-n)$
eigenvalues are of real part
$\geq\!\! (-1)$.
Hence,~$\smash{\#\widetilde{S}(\bbF_{\!\ell}) \geq \ell^2 + (2n-6)\ell + 1}$.

Among~these, at most
$(n-1)(\ell+1)$
points may have originated from blowing up the singular points
of~$S_{\frakl}$.
Indeed,~each time an
\mbox{$\bbF_{\!\ell}$-rational}
point is blown up, a
\mbox{$(+1)$-eigenspace}
is added to the Picard~group.
Therefore,
$$\#S_\reg(\bbF_{\!\ell}) \geq \ell^2 + (2n-6)\ell + 1 - (n-1)(\ell+1) = \ell^2 - 5\ell + 2 + n(\ell-1) \geq \ell^2 - 4\ell + 1 \, .$$
For~$\ell \geq 5$,
this is positive.

Thus,~it only remains to consider the case that
$\ell = 3$.
Then~$S$
is the closed subvariety
of~$\Pb^4_{\!\bbF_{\!3}}$,
given~by
\begin{eqnarray*}
T_0T_1               & = & T_2^2 - DT_3^2 \, ,\\
(T_0+aT_1)(T_0+bT_1) & = & T_2^2 - DT_4^2
\end{eqnarray*}
for
$D = \pm1$
and
certain~$a,b \in \bbF_{\!3}$.
Independently~of the values of
$a$
and~$b$,
$S$
has the regular
$\bbF_{\!3}$-rational
point
$(1\!:\!0\!:\!1\!:\!1\!:\!0)$
in the case that
$D=1$
and~$(1\!:\!0\!:\!0\!:\!0\!:\!1)$
in the case that
$D=-1$.
}
\eop
\end{coro}

\begin{prop}[Unramified primes]
\label{inert}
Let\/~$K$
be a number field,
$A,B,D \in \calO_K$,
and\/
$\frakl \subset \calO_K$
be a prime ideal that is unramified under the field extension\/
$K(\sqrt{D})/K$.
Consider~the
surface\/
$S := S^{(D;A,B)}$.

\begin{abc}
\item
If\/~$\#\calO_K/\frakl$
is not a power
of\/~$2$
then\/
$S(K_{\frakl}) \neq \emptyset$.
\item
Assume~that\/
$A \not\equiv B \pmod \frakl$,
that\/
$S$
is non-singular, and that\/
$S(K_\frakl) \neq \emptyset$.
Let\/~$\alpha \in \Br(S)$
be the Brauer class, described in Proposition~\ref{Brauerklasse}.a). Then~the local evaluation map\/
$\ev_{\alpha,\frakl} \colon S(K_\frakl) \to \bbQ/\bbZ$
is constantly~zero.
\end{abc}\smallskip

\noindent
{\bf Proof.}
{\em
We
put~$\ell := \#\calO_K/\frakl$.
Furthermore,~we normalize
$D$
to be a unit
in~$\calO_{K_\frakl}$.
This is possible because
$\frakl$
is~unramified.\smallskip

\noindent
a)
It~suffices to verify the existence of a regular
$\bbF_{\!\ell}$-rational
point on the reduction
$S_\frakl$
of~$S$.
For~this, we observe that
$(D \bmod \frakl\calO_{K_\frakl}) \neq 0$,
which shows that Corollary~\ref{ratpoint} applies.\smallskip

\noindent
b)
If~$\frakl$
is split then the assertion directly is Proposition~\ref{Brauerklasse}.b.ii). 
Otherwise,~let
$(t_0\!:\!t_1\!:\!t_2\!:\!t_3\!:\!t_4) \in S(K_\frakl)$
be an arbitrary~point. Normalize the coordinates such that
$t_0, \ldots, t_4 \in \calO_{K_\frakl}$
and at least one is a~unit.

We~first observe that one of 
$t_0$
and
$t_1$
must be a~unit. Indeed,~otherwise one has
$\frakl|t_0,t_1$.
According~to equation~(\ref{eq_eins}), this implies that
$\smash{\frakl|N_{K_\frakl(\sqrt{D})/K_\frakl}(t_2 + t_3\sqrt{D})}$.
Such a divisibility is possible only when
$\frakl|t_2,t_3$,
since
$\smash{K_\frakl(\sqrt{D})/K_\frakl}$
is an unramified, proper extension and
$\smash{\sqrt{D} \in K_\frakl(\sqrt{D})}$
is a~unit. But~then
$t_4$
is a unit, in contradiction to equation~(\ref{eq_zwei}).

Second,~we claim that
$t_0+At_1$
or
$t_0+Bt_1$
is a~unit.
Indeed,~since
$A \not\equiv B \pmod \frakl$,
the assumption
$\frakl|t_0+At_1, t_0+Bt_1$
implies
$\frakl|t_0,t_1$.

We~have thus shown that one of the four expressions
$(t_0+At_1)/t_0$,
$(t_0+At_1)/t_1$,
$(t_0+Bt_1)/t_0$,
and
$(t_0+Bt_1)/t_1$
is a~unit.
Write~$q$
for that~quotient. As~the local extension
$\smash{K_\frakl(\sqrt{D})/K_\frakl}$
is unramified of degree two, we see that
$(q,D)_\frakl = 1$.
Proposition~\ref{Brauerklasse}.b.i) implies the~assertion.
}
\eop
\end{prop}

If~$\frakl$
is a split prime then an even stronger statement is~true.

\begin{lem}[Split primes]
\label{split}
Let\/~$K$
be a number field,
$A,B,D \in \calO_K$,
and\/
$\frakl \subset \calO_K$
a prime ideal that is split under\/
$\smash{K(\sqrt{D})/K}$.
Consider~the
surface\/~$S := S^{(D;A,B)}$.

\begin{abc}
\item
Then\/
$S(K_{\frakl}) \neq \emptyset$.
\item
Furthermore,~if\/
$S$
is non-singular and\/
$\alpha \in \Br(S)$
is the Brauer class, described in Proposition \ref{Brauerklasse}.a), then the local evaluation map\/
$\ev_{\alpha,\frakl} \colon S(K_\frakl) \to \bbQ/\bbZ$
is constantly~zero.
\end{abc}\smallskip

\noindent
{\bf Proof.}
{\em
a)
The~assumption that
$\frakl$
is split under the field extension
$\smash{K(\sqrt{D})/K}$
is equivalent to
$\smash{\sqrt{D} \in K_\frakl}$.
Therefore, the point
$\smash{(1\!:\!0\!:\!1\!:\!\frac1{\sqrt{D}}\!:\!0)}$
is defined over
$K_\frakl$.
In~particular,
$S(K_\frakl)\neq \emptyset$.\smallskip

\noindent
b)
This is the assertion of Proposition~\ref{Brauerklasse}.b.ii).
}
\eop
\end{lem}

\begin{rem}
If
$\frakl$
is inert,
$0 \not\equiv A \equiv B \pmod \frakl$,
and
$(A/D \bmod \frakl) \in \calO_K/\frakl$
is a non-square then the assertion of~Proposition~\ref{inert}.b) is true,~too.

Indeed,~$t_0$
or~$t_1$
must be a unit by the same argument as~before. On~the other hand, the assumption
$\frakl|t_0+At_1, t_0+Bt_1$
does not lead to an immediate contradiction, but only to
$\frakl|t_2,t_4$
and
$t_0/t_1 \equiv -A \pmod \frakl$.
In~particular, both
$t_0$
and~$t_1$
must be~units. But~then equation~(\ref{eq_eins}) implies the congruence
$-At_1^2 \equiv -Dt_3^2 \pmod \frakl$.
\end{rem}

\begin{rem}[Inert primes--the case of residue characteristic~$2$]
\label{char2}
\leavevmode\\
We note that a statement analogous to Proposition~\ref{inert}.a) is true for any inert
prime~$\frakl$
under some more restrictive conditions on the coefficients
$A$
and
$B$.

For this suppose~that
$A$,
$B$,
$D \in \calO_K$
and that
$\frakl \subset \calO_K$
is a prime ideal that is inert under
$\smash{K(\sqrt{D})/K}$.
Let~$e$
be a positive integer such that
$x \equiv 1 \pmod {\frakl^e}$
is enough to imply that
$x \in K_\frakl$
is a~square. Assume~that
$\nu_\frakl(B-1) = f \geq 1$
and that
$\nu_\frakl(A)$
is an odd number such that
$\nu_\frakl(A) \geq 2f+e$.
Then~$S(K_\frakl) \neq \emptyset$.

Indeed, let us show that there exists a point
$(t_0\!:\!t_1\!:\!t_2\!:\!t_3\!:\!t_4) \in S(K_\frakl)$
such that
$t_3=t_4$
and~$t_1 \neq 0$.
This~leads to the equation
$(T_0+AT_1)(T_0+BT_1) = T_0T_1$,
or
$$T_0^2 + (A+B-1)T_0T_1 + ABT_1^2 = 0 \, .$$
The~discriminant of this binary quadric is
$(A+B-1)^2 - 4AB = (B\!-\!1)^2 + A(A-2B\!-\!2)$,
which is a square
in~$K_\frakl$
by virtue of our~assumptions. Thus,~there are two solutions
in~$K_\frakl$
for
$T_0/T_1$
and their product is
$AB$,
which is of odd~valuation. We~may therefore choose a solution
$t_0/t_1$
such that
$\nu_\frakl(t_0/t_1)$
is~even. This~is enough to imply that
$(t_0+At_1)(t_0+Bt_1) = t_0t_1$
is a norm
from~$\smash{K_\frakl(\sqrt{D})}$.
\end{rem}

\begin{rem}[Archimedean primes]
\label{Archimedean}
\begin{iii}
\item
Let~$\sigma\colon K \hookrightarrow \bbR$
be a real~prime. Then,~for
$A$,
$B \in K$
ar\-bitrary and
$D \in K$
non-zero, one has
$S_\sigma(\bbR) \neq \emptyset$.

Indeed,~we can put
$t_1 := 1$
and choose
$t_0 \in \bbR$
such that
$t_0$,
$t_0 + \sigma(A)$,
and
$t_0 + \sigma(B)$
are~positive. Then
$C := t_0 > 0$
and
$C' := (t_0 + \sigma(A))(t_0 + \sigma(B)) > 0$
and we have to show that the system of equations
\begin{eqnarray*}
T_2^2 - \sigma(D) T_3^2 & = & C \\
T_2^2 - \sigma(D) T_4^2 & = & C'
\end{eqnarray*}
is solvable in~$\bbR$.
For this one may choose
$t_2$
such that
$t_2^2 \geq \max(C,C')$
if
$\sigma(D) > 0$
and such that
$t_2^2 \leq \min(C,C')$,
otherwise. In~both cases it is clear that there exist real numbers
$t_3$
and
$t_4$
such that the resulting point is contained
in~$S_\sigma(\bbR)$.

Moreover~if
$\sigma(D) > 0$
then the local evaluation map
$\ev_{\alpha,\sigma}\colon S(K_\sigma) \to \frac12\bbZ/\bbZ$
is constantly~zero. Indeed,~then one has
$(q,D)_\sigma = 1$
for every
$q \in K_\sigma \cong \bbR$,
different from~zero.
\item
For~$\tau\colon K \hookrightarrow \bbC$
a complex prime and
$A$,
$B$,
and
$D \in K$
arbitrary, we clearly have that
$S(K_\tau) \neq \emptyset$.
Furthermore,~$(q,D)_\tau = 1$
for every non-zero
$q \in K_\tau \cong \bbC$.
\end{iii}
\end{rem}

\section{Ramification--Reduction to the union of four planes}

The goal of this section is to study the evaluation of the Brauer class at ramified primes
$\frakl$.
Under~certain congruence conditions on the parameters
$A$
and~$B$
we deduce that the evaluation map is constant on the
$K_\frakl$-rational
points on
$S$,
and we determine its value depending on
$A$
and~$B$.

\begin{prop}[Ramified primes in residue characteristic $\neq \!2$]
\label{ramified}
\leavevmode\\
Let\/~$K$
be a number field,
$A,B,D \in \calO_K$,
and\/
$\frakl \subset \calO_K$
a prime ideal such that\/
$\#\calO_K/\frakl$
is not a power
of\/~$2$
and that is ramified under the field extension\/
$\smash{K(\sqrt{D})/K}$.
Suppose~that\/
\mbox{$\overline{A} := (A \bmod \frakl)$}
$\in \calO_K/\frakl$
is a square, different from\/
$0$
and\/
$(-1)$,
that\/
$\smash{\overline{A}^2 + \overline{A} + 1 \neq 0}$,
and~that
$$\textstyle B \equiv -\frac{A}{A+1} \pmod \frakl \, .$$
Consider~the
surface\/~$S := S^{(D;A,B)}$.

\begin{abc}
\item
Then\/
$S(K_{\frakl}) \neq \emptyset$.
\item
Assume~that\/
$S$
is non-singular and let\/~$\alpha \in \Br(S)$
be the Brauer class, described in Proposition \ref{Brauerklasse}.a).
\begin{iii}
\item
If\/
$\overline{A}+1 \in \calO_K/\frakl$
is a square then the local evaluation map\/
$\ev_{\alpha,\frakl}\colon S(K_\frakl) \to \bbQ/\bbZ$
is constantly~zero.
\item
If\/
$\overline{A}+1 \in \calO_K\!/\frakl$
is a non-square then the local evaluation map\/
$\ev_{\alpha,\frakl}\colon S(K_\frakl) \to \bbQ/\bbZ$
is constant~of
value~$\frac12$.
\end{iii}
\end{abc}\smallskip

\noindent
{\bf Proof.}
{\em
First~of all, we note that
$\nu_\frakl(D)$
is~odd.
Indeed,~assume the contrary. We~may then normalize
$D$
to be a unit and write
$K_\frakl^\n$
for the unramified quadratic extension
of~$K_\frakl$.
Then~$(D \bmod \frakl\calO_{K_\frakl^\n})$
is a square and, since
$\calO_{K_\frakl^\n}/\frakl\calO_{K_\frakl^\n}$
is a field of characteristic different
from~$2$,
Hensel's Lemma ensures that
$D$
is a square in
$K_\frakl^\n$.
I.e.,~$\smash{K_\frakl(\sqrt{D}) \subseteq K_\frakl^\n}$,
a~contradiction.

Let~us normalize
$D$
such that
$\nu_\frakl(D) = 1$.
Then~the reduction
$S_\frakl$
of~$S$
is given by the equations
\begin{eqnarray}
T_0T_1
& = & T_2^2 \, , \label{eq_drei} \\
\textstyle
\smash{(T_0+\overline{A}T_1)(T_0-\frac{\overline{A}}{\overline{A}+1}T_1)}
& = & T_2^2 \, , \label{eq_vier}
\end{eqnarray}
which geometrically define a cone over a cone over four points
in~$\Pb^2$.\smallskip

\noindent
a)
We~write
$\ell := \#\calO_K/\frakl$.
It~suffices to verify the existence of a regular
$\bbF_{\!\ell}$-rational
point
on~$S_\frakl$.
For~this, it is clearly enough to show that one of the four points
in~$\Pb^2$,
defined by the equations (\ref{eq_drei}) and~(\ref{eq_vier}), is simple and defined
over~$\bbF_{\!\ell}$.

Equating the two terms on the left hand side, one finds the equation
$$\textstyle T_0^2 + \frac{\overline{A}^2-\overline{A}-1}{\overline{A}+1}T_0T_1 - \frac{\overline{A}^2}{\overline{A}+1}T_1^2 = 0 \, ,$$
which obviously has the two solutions
$T_0/T_1 = 1$
and
$\smash{T_0/T_1 = - \frac{\overline{A}^2}{\overline{A}+1}}$.
By~virtue of our assumptions, both are
\mbox{$\bbF_{\!\ell}$-rational}
points
in~$\Pb^1$,
different from
$0$
and
$\infty$.
They~are different from each other, since
$\smash{\overline{A}^2 + \overline{A} + 1 \neq 0}$.

Consequently,~the four points defined by the equations (\ref{eq_drei}) and~(\ref{eq_vier}) are all~simple. The~two points corresponding to
$(t_0\!:\!t_1) = 1$
are defined
over~$\bbF_{\!\ell}$.
The~two others are defined
over~$\bbF_{\!\ell}$
if and only if
$(-\overline{A}-1) \in \bbF_{\!\ell}$
is a~square.\smallskip

\noindent
b)
Let
$(t_0\!:\!t_1\!:\!t_2\!:\!t_3\!:\!t_4) \in S(K_\frakl)$
be any point. We~normalize the coordinates such that
$t_0,\ldots,t_4 \in \calO_{K_l}$
and at least one of them is a~unit.
Then
$\frakl$
cannot divide both
$t_0$
and
$t_1$.
Indeed,~this would imply
$\frakl^2 | t_2^2 - Dt_3^2$
and
$\frakl^2 | t_2^2 - Dt_4^2$
and,
as~$\nu_\frakl(D) = 1$,
this is possible only for
$\frakl | t_2, t_3, t_4$.

Therefore,~$((t_0 \!+\! At_1)/t_1 \bmod \frakl) = \overline{A} \!+\! (t_0/t_1 \bmod \frakl)$
is either equal to
$\smash{(\overline{A} \!+\! 1)}$
or to
$\smash{\overline{A} \!-\! \frac{\overline{A}^2}{\overline{A}+1} = \frac{\overline{A}}{\overline{A}+1}}$.
Both~terms are squares
in~$\bbF_{\!\ell}$
under the assumptions of b.i), while, under the assumptions of~b.ii), both are non-squares.

As~a unit
in~$\calO_{K_\frakl}$
is a norm from the ramified
extension~$\smash{K_\frakl(\sqrt{D})}$
if and only if its residue
modulo~$\frakl$
is a square, for
$q := (t_0+At_1)/t_1$,
we find that
$(q,D)_\frakl = 1$
in case~i) and
$(q,D)_\frakl = -1$
in case~ii). Proposition~\ref{Brauerklasse}.b.ii) implies the~assertion.
}
\eop
\end{prop}

\section{The main result}

We are now in the position to formulate sufficient conditions on
$A,B,D$,
under which the corresponding surface
$S^{(D;A,B)}$
violates the Hasse principle.

\begin{theo}
\label{Data_Hasse_Gegen}
Let\/~$D \in K$
be non-zero and\/
$(D) = (\frakq_1^{k_1}\cdot\ldots\cdot\frakq_l^{k_l})^2\, \frakp_1\cdot\ldots\cdot\frakp_k$
its decomposition into prime~ideals. Suppose~that the primes\/
$\frakp_1, \ldots, \frakp_k$
are~distinct.

\begin{abc}
\item
Assume~that
\begin{iii}
\item
$k \geq 1$,
\item
the quadratic extension\/
$\smash{K(\sqrt{D})/K}$
is unramified at all primes
of\/~$K$
lying over the rational
prime\/~$2$,
\item
for~every real prime\/
$\sigma\colon K \hookrightarrow \bbR$,
one has\/
$\sigma(D) > 0$.
\end{iii}
\item
For~every
prime\/~$\frakl$
of\/~$K$
that lies over the rational
prime\/~$2$
and is inert under\/
$\smash{K(\sqrt{D})/K}$,
assume~that
$$\nu_\frakl(B-1) = f_\frakl \geq 1\,, \quad{\it that\/~} \nu_\frakl(A) {\it ~is~odd \,, ~} \quad {\it and~that\/~} \nu_\frakl(A) \geq 2f_\frakl+e_\frakl\,,$$
for\/~$e_\frakl$
a positive integer such that\/
$x \equiv 1 \pmod {\frakl^{e_\frakl}}$
is enough to ensure that\/
$x \in K_\frakl$
is a~square.
\item
For~every\/
$i = 1,\ldots,k$,
assume~that
\begin{iii}
\item
$(A \bmod \frakp_i) \in \calO_K/\frakp_i$
is a square, different from\/
$0$,
$(-1)$,
and the primitive third roots of~unity.
If\/~$\#\calO_K/\frakp_i$
is a power
of\/~$3$
then assume\/
$(A \bmod \frakp_i) \neq 1$,~too.
\item
$B \equiv -\frac{A}{A+1} \pmod {\frakp_i}$.
\item
$1+(A \bmod \frakp_i) \in \calO_K/\frakp_i$
is a non-square for\/
$i=1,\ldots,b$,
for an odd
integer\/~$b$,
and a square
for\/
$i=b+1,\ldots,k$.
\end{iii}
\item
Finally,~assume that\/
$(A-B)$
is a product of only split~primes.
\end{abc}\smallskip

\noindent
Then\/
$S^{(D;A,B)}(\bbA_K) \neq \emptyset$.
However,~if\/
$S^{(D;A,B)}$
is non-singular then\/
$S^{(D;A,B)}(K) = \emptyset$.
\end{theo}

\begin{rem}
Without any change, one may assume that the 
$\frakq_1, \ldots, \frakq_l$
are distinct,~too. Note,~however, that we do not suppose
$\{\frakp_1, \ldots, \frakp_k\}$
and
$\{\frakq_1, \ldots, \frakq_l\}$
to be disjoint.
\end{rem}

\noindent
{\bf Proof of Theorem~\ref{Data_Hasse_Gegen}.}
By~a.i),
$D$
is not a square
in~$K$,
hence
$K(\sqrt{D})/K$
is a proper quadratic field~extension. It~is clearly ramified at
$\frakp_1, \ldots, \frakp_k$.
According~to a.ii), these are the only ramified~primes.
In~view of assumption~b),
$S(\bbA_\bbQ) \neq \emptyset$
follows from Proposition~\ref{inert}.a) and Proposition~\ref{ramified}.a), together with Lemma \ref{split}.a), Remark \ref{char2}, and Remark~\ref{Archimedean}.

On~the other hand, let
$\alpha \in \Br(S)$
be the Brauer class, described in Proposition~\ref{Brauerklasse}.a). Then,~in view of assumptions~d), c) and~a.iii), Proposition~\ref{inert}.b) and Proposition~\ref{ramified}.b), together with Lemma  \ref{split}.b) and Remark~\ref{Archimedean}, show that the local evaluation map
$\ev_{\alpha,\frakl}$
is constant of
value~$\frac12$
for
$\frakl = \frakp_1,\ldots,\frakp_b$
and constantly zero for all others. Proposition~\ref{BM} proves that
$S$
is a counterexample to the Hasse~principle.
\eop

\begin{ex}
Let~$S$
be the surface
in~$\Pb^4_\bbQ$,
given~by
\begin{eqnarray*}
               T_0T_1 & = & T_2^2 - 17T_3^2 \, ,\\
(T_0+9T_1)(T_0+11T_1) & = & T_2^2 - 17T_4^2 \, .
\end{eqnarray*}
Then
$S(\bbA_\bbQ) \neq \emptyset$
but~$S(\bbQ) = \emptyset$.\medskip

\noindent
{\bf Proof.}
We~have
$K=\bbQ$
and
$D=17$.
Furthermore,
$A = 9$
and
$B = 11$
such that Proposition~\ref{nonsing} ensures that
$S = S^{(D;A,B)}$
is non-singular.

The~extension
$\smash{L := \bbQ(\sqrt{17})/\bbQ}$
is real-quadratic,
i.e.~$D > 0$,
and ramified only
at~$17$.
Under~$\smash{\bbQ(\sqrt{17})/\bbQ}$,
the prime
$2$
is split, which completes the verification of~a) and shows that b) is fulfilled~trivially.

For~c), note that
$17 \not\equiv 1 \pmod 3$,
such that there are no nontrivial third roots of unity
in~$\bbF_{\!17}$.
Furthermore,
$9 \neq 0,(-1)$
is a square
modulo~$17$,
but~$10$
is not,
and~$\smash{11 \equiv -\frac9{10} \pmod {17}}$.
Finally, for~d), note that
$(A-B) = (-2) = (2)$
is a prime that is split
in~$\smash{\bbQ(\sqrt{17})}$.
\eop
\end{ex}

\begin{rem}
The assumption
on~$S$
to be non-singular may be removed from Theorem~\ref{Data_Hasse_Gegen}. Indeed, the elementary argument described at the very end of section~\ref{sec_Brauer} works in the singular case,~too.
\end{rem}

The goal of the next lemma is to construct discriminants
$D \in K$,
for which we will later be able to construct counterexamples to the Hasse principle, via the previous theorem.

\begin{lem}
\label{D_exists}
Let\/~$K$
be an arbitrary number field and\/
$\frakp, \frakr_1, \ldots, \frakr_n$
be distinct prime ideals such that\/
$\calO_K/\frakp$
and\/~$\calO_K/\frakr_i$
are of characteristics different
from\/~$2$.
Then~there exists some\/
$D \in K$
such~that

\begin{iii}
\item
the prime\/
$\frakp$
is ramified
in\/~$\smash{K(\sqrt{D})}$,
\item
all primes lying over the rational
prime\/~$2$
are split
in\/~$\smash{K(\sqrt{D})}$.
\item
For~every real prime\/
$\sigma\colon K \hookrightarrow \bbR$,
one has\/
$\sigma(D) > 0$.
\item
The~primes\/
$\frakr_i$
are unramified
in\/~$\smash{K(\sqrt{D})}$.
\end{iii}

\noindent
In~particular, assumptions~a) and~b) of Theorem~\ref{Data_Hasse_Gegen} are~fulfilled.\medskip

\noindent
{\bf Proof.}
{\em
Let
$\frakl_1,\ldots,\frakl_m$
be the primes
of~$K$
that lie over the rational
prime~$2$.
We~impose the congruence conditions
$D \equiv 1 \pmod {\frakl_1^{e_1}}, \ldots, D \equiv 1 \pmod {\frakl_m^{e_m}}$,
for~$e_1, \ldots, e_m$
large enough that this implies that
$D$
is a square
in~$K_{\frakl_1}, \ldots, K_{\frakl_m}$.

Furthermore,~the assumptions imply that
$\frakp, \frakr_1, \ldots, \frakr_n$
are different from
$\frakl_1,\ldots,\frakl_m$.
We~impose, in addition, the conditions
$D \in \frakp \!\setminus\! \frakp^2$
and
$D \not\in \frakr_1, \ldots, \frakr_n$.

According to the Chinese remainder theorem, these conditions have a simultaneous
solution~$D'$.
Put
$D := D' + k\cdot\#\calO_K/\frakl_1^{e_1}\ldots\frakl_m^{e_m}\frakp^2\frakr_1\ldots\frakr_n$,
for~$k$
an integer that is sufficiently large to ensure
$\sigma(D) > 0$
for every real prime
$\sigma\colon K \hookrightarrow \bbR$.
Then~assertion~iii) is~true.
Furthermore,~the congruences
$D \equiv 1 \pmod {\frakl_i^{e_i}}$
imply~ii), while
$D \in \frakp \!\setminus\! \frakp^2$
yields assertion~i) and 
$D \not\in \frakr_1, \ldots, \frakr_n$
ensures that~iv) is~true.
}
\eop
\end{lem}

Before we come to the next main theorem of this section, we need to formulate two technical lemmata. 

\begin{lem}
\label{Chebotarev}
Let\/~$K$
be a number field,
$I \subset \calO_K$
an ideal,
and\/~$x \in \calO_K$
an element relatively prime
to\/~$I$.\smallskip

\noindent
Then~there exists an infinite sequence of pairwise non-associated elements\/
$y_i \in \calO_K$
such that, for
each\/~$i \in \bbN$,
one has that\/
$(y_i)$
is a prime ideal and\/
$y_i \equiv x \pmod I$.\medskip

\noindent
{\bf Proof.}
{\em
It~is well known that there exist infinitely many prime ideals\/
$\frakr_i \subset \calO_K$
with the property~below.\smallskip

There~exist some\/
$u_i, v_i \in\calO_K$,
$u_i \equiv v_i \equiv 1 \pmod I$
such~that
$$\frakr_i \!\cdot\! (u_i) = (x) \!\cdot\! (v_i) \, .\medskip$$

\noindent
Indeed,~the invertible ideals
in~$K$
modulo the principal ideals generated by elements from the residue class
$(1 \bmod I)$
form an abelian group that is canonically isomorphic to the ray class group
$Cl_K^I \cong C_K/C_K^I$
of~$K$
\cite[Chapter~VI, Proposition~1.9]{Ne}. Thus,~the claim follows from the Chebotarev density theorem applied to the ray class field
$K^I\!/\!K$,
which has the Galois group
$\Gal(K^I\!/\!K) \cong Cl_K^I$.

Take~one of these prime~ideals.
Then~$\frakr_i \!\cdot\! (u_i) = (x) \!\cdot\! (v_i) = (xv_i)$.
As~$\frakr_i \subset \calO_K$,
this shows that
$xv_i$
is divisible
by~$u_i$.
Put~$y_i := xv_i/u_i$.
Then~$(y_i) = \frakr_i$.
Further,
\mbox{$y_i \equiv x \pmod I$}.
}
\eop
\end{lem}

\begin{lem}
\label{sq_nonsq_ex}
Let\/~$\bbF_{\!q}$
be a finite field of characteristic\/
$\neq \!2$
having\/
$> \!25$~elements.
Then~there exist elements\/
$a_{00}$,
$a_{01}$,
$a_{10}$,
and
$a_{11} \in \bbF_{\!q}$,
different from\/
$0$,
$(-1)$,
$(-2)$
and such that\/
$a_{ij}^2 + a_{ij} + 1 \neq 0$,
that fulfill the conditions~below.

\begin{iii}
\item
$a_{00}$,
$(a_{00}+1)$,
and\/
$(a_{00}+2)$
are squares
in\/~$\bbF_{\!q}$.
\item
$a_{01}$
and\/
$(a_{01}+1)$
are squares
in\/~$\bbF_{\!q}$,
but\/
$(a_{01}+2)$
is~not.
\item
$a_{10}$
and\/
$(a_{10}+2)$
are squares
in\/~$\bbF_{\!q}$,
but\/
$(a_{10}+1)$
is~not.
\item
$a_{11}$
is a square
in\/~$\bbF_{\!q}$,
but\/
$(a_{11}+1)$
and\/
$(a_{11}+2)$
are~not.
\end{iii}\smallskip

\noindent
{\bf Proof.}
{\em
Let~$\smash{C_1 \in \bbF_{\!q}^*}$
be a square in the cases i) and ii), and a non-square,~otherwise. Similarly,~let
$\smash{C_2 \in \bbF_{\!q}^*}$
be a square in the cases i) and iii), and a non-square,~otherwise. The~problem then translates into finding an
\mbox{$\bbF_{\!q}$-rational}
point on the
curve~$E$,
given
in~$\Pb^3$~by
\begin{eqnarray*}
U_1^2 + \phantom{1}U_0^2 & = & C_1U_2^2 \,, \\
U_1^2 +           2U_0^2 & = & C_2U_3^2 \,,
\end{eqnarray*}
such that
$U_i \neq 0$
for
$i = 0,\ldots,3$
and
$(U_1/U_0)^4 + (U_1/U_0)^2 + 1 \neq 0$.
Note~that the conditions
$U_2 \neq 0$
and
$U_3 \neq 0$
imply that
$\smash{\big(\frac{U_1}{U_0}\big)^2 \neq -1,-2}$.

Since~the characteristic of the base field is different from two, a direct calculation shows that
$E$~is
non-singular, i.e.~a smooth curve of
genus~$1$.
The~extra conditions define an open subscheme
$\smash{\widetilde{E} \subset E}$
that excludes not more than 32~points. Thus,~Hasse's bound yields
$\smash{\#\widetilde{E}(\bbF_{\!q}) \geq q - 2\sqrt{\mathstrut q} - 31}$.
This~is positive
for~$q > 44$.

An~experiment shows that
$\bbF_{\!3}$,
$\bbF_{\!5}$,
$\bbF_{\!7}$,
$\bbF_{\!9}$,
$\bbF_{\!13}$,
$\bbF_{\!17}$,
and
$\bbF_{\!25}$
are the only fields in
characteristic~$\neq \!2$,
for which the assertion is~false.
}
\eop
\end{lem}

The following theorem provides us with Hasse counterexamples in the family $S^{(D;A,B)}$
for suitable
discriminants~$D$.
For us, the important feature is that one may choose the parameters
$A$
and~$B$
to lie in (almost) arbitrary congruence classes modulo some prime ideal
$\frakl \subset \calO_K$,
unramified
in
$K(\sqrt{D})$,
provided only that
$A \not\equiv B \pmod \frakl$. 

\begin{theo}
\label{Hasse_congr_ex}
Let\/~$K$
be an arbitrary number field and\/
$D \in K$
a non-zero~element. Write\/
$(D) = (\frakq_1^{k_1}\cdot\ldots\cdot\frakq_l^{k_l})^2 \frakp_1\cdot\ldots\cdot\frakp_k$
for its decomposition into prime ideals, the\/
$\frakp_i$
being~distinct. Assume~that

\begin{iii}
\item
$k \geq 1$,
\item
all primes lying over the rational
prime\/~$2$
are split
in\/~$\smash{K(\sqrt{D})}$,
\item
for~every real prime\/
$\sigma\colon K \hookrightarrow \bbR$,
one has\/
$\sigma(D) > 0$,
\item
all primes with residue field\/
$\bbF_{\!3}$
are unramified
in\/~$\smash{K(\sqrt{D})}$.
\end{iii}

\noindent
Suppose~further that among the primes\/
$\frakp$
of\/~$K$
that are ramified
in\/~$\smash{K(\sqrt{D})}$,
there is one such that\/
$\#\calO_K/\frakp > 25$.\smallskip

\noindent
Then,~for every prime\/
$\frakl \subset \calO_K$,
unramified
in\/~$\smash{K(\sqrt{D})}$,
and all\/
$a,b \in \calO_K/\frakl$
such that\/
$a \neq b$,
there exist
$A,B \in \calO_K$
such that\/
$(A \bmod \frakl) = a$,
$(B \bmod \frakl) = b$, and $S^{(D;A,B)}(\bbA_K) \neq \emptyset$,
but\/
$S^{(D;A,B)}(K) = \emptyset$.\medskip

\noindent
{\bf Proof.}
{\em
{\em First step.}
Construction of
$A$
and~$B$.

\noindent
Let~$M \in \{1,\ldots,k\}$
be such that
$\#\calO_K/\frakp_M > 25$.
Besides
\begin{equation}
\label{cong1}
(A \bmod \frakl) = a \qquad {\rm and} \qquad (B \bmod \frakl) = b \, ,
\end{equation}
we will impose further congruence~conditions on
$A$
and~$B$.
For~each
$i \neq M$,
we choose a square
$a_i \in \calO_K/\frakp_i$
such that
$a_i \neq 0, (-1), (-2)$
and
$a_i^2 + a_i + 1 \neq 0$.
This~is possible since
$\calO_K/\frakp_i$
is of
char\-ac\-ter\-is\-tic~$\neq \!2$
and~$\#\calO_K/\frakp_i > 3$.
We~require
\begin{equation}
\label{cong2}
(A \bmod \frakp_i) = a_i \qquad {\rm and} \qquad (B \bmod \frakp_i) = -\frac{a_i}{a_i+1} \, .
\end{equation}
Finally,~we choose a square
$a_M \in \calO_K/\frakp_M$
such that
$a_M \neq 0, (-1), (-2)$
and
$a_M^2 + a_M + 1 \neq 0$,
satisfying the additional conditions~below.

\begin{iii}
\item[$\bullet$ ]
If,~among the elements
$a_1+1, \ldots, a_{M-1}+1, a_{M+1}+1, \ldots, a_k+1$,
there are an odd number of non-squares then
$a_M+1$
is a~square.
Otherwise,~$a_M+1$
is a non-square.
\item[$\bullet$ ]
If,~among the elements
$a_1+2, \ldots, a_{M-1}+2, a_{M+1}+2, \ldots, a_k+2$,
there are an odd number of non-squares then
$a_M+2$
is a~square.
Otherwise,~$a_M+2$
is a non-square.
\end{iii}
Lemma~\ref{sq_nonsq_ex} guarantees that such an element
$a_M \in \calO_K/\frakp_M$~exists.
We~impose the final congruence~condition
\begin{equation}
\label{cong3}
(A \bmod \frakp_M) = a_M \qquad {\rm and} \qquad (B \bmod \frakp_M) = -\frac{a_M}{a_M+1} \, .
\end{equation}\smallskip

According~to the Chinese remainder theorem, one may choose an algebraic integer
$B \in \calO_K$
such that the conditions on the right hand sides of (\ref{cong1}), (\ref{cong2}), and (\ref{cong3}) are~fulfilled. Then,~by Lemma~\ref{Chebotarev}, there exist infinitely many non-associated elements
$y_i \in \calO_K$
such that
$(y_i)$
is a prime ideal and
$(y_i+B,B)$
a simultaneous solution of the system of congruences~(\ref{cong1},\,\ref{cong2},\,\ref{cong3}).

We~choose some
$i \in \bbN$
such that
$\frakr := (y_i)$
is of residue characteristic different
from~$2$,
that
$\frakr \neq \frakp_1, \ldots, \frakp_k, \frakq_1, \ldots, \frakq_l$,
and such that
$A^2-2AB+B^2-2A-2B+1 \neq 0$
for
$A := y_i+B$.
Note~that
$\frakr \neq \frakp_1, \ldots, \frakp_k, \frakq_1, \ldots, \frakq_l$
is equivalent to
$\frakr \not\ni D$.\smallskip

\noindent
{\em Second step.}
The surface
$S := S^{(D;A,B)}$
is a counterexample to the Hasse~principle.

\noindent
To~show this, let us use Theorem~\ref{Data_Hasse_Gegen}.
Our~assumptions
on~$D$
imply that assumptions~a) and~b) of Theorem~\ref{Data_Hasse_Gegen} are~fulfilled. Assumption~c) is satisfied, too, by consequence of the construction of the\
elements~$a_i$.
Observe,~in particular, that among the elements
$a_1+1, \ldots, a_k+1$,
there are an odd number of~non-squares.
Furthermore,~$S$
is non-singular.

It~therefore remains to check assumption~d). The~only
prime~$\frakp \subset \calO_K$,
for which
$A \equiv B \pmod \frakp$,
is
$\frakp = \frakr \;\,(=\!(A\!-\!B))$.
We~have to show that
$\frakr$
is split
under~$\smash{K(\sqrt{D})/K}$.\smallskip

For~this, we observe that,
for~$i = 1,\ldots,k$,
$$\textstyle A-B \equiv A + \frac{A}{A+1} = A\frac{A+2}{A+1} \pmod {\frakp_i} \, .$$
As~$A$
is a square
modulo~$\frakp_i$,
this shows
$$\prod_{i=1}^k (A-B,D)_{\frakp_i} = \prod_{i=1}^k (A+2,D)_{\frakp_i} \Big/ \prod_{i=1}^k (A+1,D)_{\frakp_i} \, .$$
Here,~by our construction, both
$1 + (A \bmod \frakp_i)$
and
$2 + (A \bmod \frakp_i)$
are non-squares, an odd number of~times. Consequently,
$$\prod_{i=1}^k (A-B,D)_{\frakp_i} = 1 \, .$$\smallskip

On~the other hand,
$D$~is
a square
in~$K_{\frakl_i}$
for~$l_i$
the primes of residue
characteristic~$2$
and for every real~prime, by assumption iii).
Thus,~$(A-B,D)_\frakl = 1$
unless
$\frakl$
divides either
$(A-B)$
or~$D$.
I.e.~for
$\frakl \neq \frakr, \frakp_1, \ldots, \frakp_k, \frakq_1, \ldots, \frakq_l$.
Moreover,~$(A-B,D)_\frakq = 1$
for~$\frakq \in \{\frakq_1, \ldots, \frakq_l\} \!\setminus\! \{\frakp_1, \ldots, \frakp_k\}$
since both arguments of the Hilbert symbol are of even
\mbox{$\frakq$-adic}
valuation.
The~Hilbert reciprocity law~\cite[Chapter~VI, Theorem~8.1]{Ne} therefore reveals the fact that
$$(A-B,D)_\frakr \cdot \prod_{i=1}^k (A-B,D)_{\frakp_i} = 1 \, .$$

Altogether, this implies
$(A-B,D)_\frakr = 1$.
Consequently,~the prime ideal
$\frakr$
splits
in~$\smash{K(\sqrt{D})}$.
}
\eop
\end{theo}

\begin{subl}
\label{01infty}
The~rational map
$\kappa\colon \Ab^{\!2}/S_2 \ratarrow \calU\!/(S_5 \times \PGL_2)$,
given on points~by
$$\overline{(a_1,a_2)} \mapsto \overline{(a_1,a_2,0,-1,\infty)} \, ,$$
is~dominant.\medskip

\noindent
{\bf Proof.}
{\em
It~suffices to prove that the rational map
$\widetilde\kappa\colon \Ab^{\!2} \ratarrow \calU\!/(S_5 \times \PGL_2)$,
given by
$\smash{(a_1,a_2) \mapsto \overline{(a_1,a_2,0,-1,\infty)}}$
is~dominant. For~this, recall that dominance may be tested after base extension to the algebraic~closure. Moreover,~it is well known that three distinct points
on~$\smash{\Pb^1_{\overline{K}}}$
may be sent to
$0$,~$(-1)$,
and~$\infty$
under the operation
of~$\smash{\PGL_2(\overline{K})}$.
}
\eop
\end{subl}

\begin{lem}
\label{inv_dom}
Let\/~$K$
be any field of
characteristic\/~$\neq \!2$
and\/
$0 \neq D \in K\!$.
Let\/
$\pi\colon \calS \to U$
be the family of degree four del Pezzo surfaces over an open subscheme\/
$U \subset \Ab^{\!2}_K$,
given by
\begin{eqnarray*}
                  T_0T_1 & = & T_2^2 - DT_3^2 \, ,\\
(T_0+a_1T_1)(T_0+a_2T_1) & = & T_2^2 - DT_4^2 \, .
\end{eqnarray*}
I.e.,~the fiber of\/
$\pi$
over\/
$(a_1,a_2)$
is exactly the surface\/
$S^{(D;a_1,a_2)}$. 
Then~the invariant map
$$I_\pi\colon U \longrightarrow\Pb(1,2,3) $$
associated
with\/~$\pi$
is~dominant.\medskip

\noindent
{\bf Proof.}
{\em
As~dominance may be tested after base extension to the algebraic closure, let us assume that the base
field~$K$
is algebraically~closed. Write
\begin{eqnarray*}
Q_1(a_1, a_2; T_0,\ldots,T_4) & := & T_0T_1 - (T_2^2 - DT_3^2) \hspace{3.5cm} {\rm and} \\
Q_2(a_1, a_2; T_0,\ldots,T_4) & := & (T_0+a_1T_1)(T_0+a_2T_1) - (T_2^2 - DT_4^2), \,
\end{eqnarray*}
and consider the pencil
$\smash{(uQ_1+vQ_2)_{(u:v) \in \Pb^1}}$
of quadrics, parametrized by
$(a_1, a_2) \in \Ab^{\!2}(K)$.

We~see that, independently of the values of the parameters, degenerate quadrics occur for
$(u:v) = 0$,
$\infty$,
and~$(-1)$.
The~two other degenerate quadrics appear for
$(u:v)$
the zeroes of the determinant
$$
\left|
\begin{array}{cc}
1             & (a_1+a_2+t)/2 \\
(a_1+a_2+t)/2 & a_1a_2
\end{array}
\right|
= \textstyle -\frac14 [t^2 + 2(a_1+a_2)t + (a_1-a_2)^2] \, .
$$
Thus,~$I_\pi$
is the composition of the rational map
$\rho\colon \Ab^{\!2} \supset U \ratarrow \Ab^{\!2}/S_2$,
sending
$(a_1,a_2)$
to the pair of roots
of~$t^2 + 2(a_1+a_2)t + (a_1-a_2)^2$,
followed by the rational map
$\kappa\colon \Ab^{\!2}/S_2 \ratarrow \calU\!/(S_5 \times \PGL_2)$,
studied in Sublemma~\ref{01infty}, and the open embedding
$\iota\colon \calU\!/(S_5 \times \PGL_2) \hookrightarrow \Pb(1,2,3)$,
defined by the fundamental~invariants. It~remains to prove that
$\rho\colon U \ratarrow \Ab^{\!2}/S_2$
is~dominant.

For~this, as coordinates on
$\Ab^{\!2}/S_2$
one may choose the sum and the product of the coordinates
on~$\Ab^{\!2}$.
Indeed,~these generate the field of
\mbox{$S_2$-invariant}
functions
on~$\Ab^{\!2}$.
Thus,~we actually claim that the map
$\Ab^{\!2} \to \Ab^{\!2}$, given by
$(a_1,a_2) \mapsto (-2(a_1+a_2), (a_1-a_2)^2)$
is dominant, which is~obvious.
}
\eop
\end{lem}

We~are now, finally, in the position to prove that the set of counterexamples to the Hasse principle is Zariski dense in the moduli scheme of del Pezzo surfaces of degree four. For~this, we will consider the family
$S^{(D;A,B)}$
for some fixed
discriminant~$D$
and use Theorem~\ref{Hasse_congr_ex}. It~will turn out to provide us with enough counterexamples to rule out the possibility that all pairs
$(A,B)$
leading to counterexamples might be contained in some finite union of~curves. Moreover, Lemma~\ref{inv_dom} shows that, for some fixed
$D \neq 0$,
the family
$S^{(D;A,B)}$
already dominates the moduli scheme of all degree four del Pezzo~surfaces.

\begin{theo}
Let\/~$K$
be any number field,
$U_\reg \subset \Gr(2,15)_K$
the open subset of the Gra\ss mann scheme that parametrizes degree four del Pezzo surfaces, and\/
$\calHC_K \subset U_\reg(K)$
be the set of all degree four del Pezzo surfaces
over\/~$K$
that are counterexamples to the Hasse~principle.\smallskip

\noindent
Then~the image of\/
$\calHC_K$
under the invariant map
$$\smash{I\colon U_\reg \longrightarrow \Pb(1,2,3)_K}$$
is Zariski~dense.\medskip

\noindent
{\bf Proof.}
{\em
According~to Lemma~\ref{D_exists}, there exists an algebraic integer
$D \in \calO_K$
fulfilling the assumptions of Theorem~\ref{Hasse_congr_ex}. Assume~that the image
of~$I$
would not be Zariski~dense. By~Lemma~\ref{inv_dom}, this implies that there exists a (possibly reducible) curve
$C \subset \Ab^{\!2}$
of certain
degree~$d$
such that, for all surfaces of the~form
\begin{eqnarray*}
T_0T_1               & = & T_2^2 - DT_3^2 \, , \\
(T_0+AT_1)(T_0+BT_1) & = & T_2^2 - DT_4^2
\end{eqnarray*}
that violate the Hasse principle, one
has~$(A,B) \in C(K)$.

On~the other hand, let
$\frakl \subset \calO_K$
be an unramified prime and put
$\ell := \#\calO_K/\frakl$.
Then,~by Theorem~\ref{Hasse_congr_ex}, we know counterexamples to the Hasse principle having
$\ell(\ell-1)$
distinct reductions
modulo~$\frakl$.
But~an affine plane curve of
degree~$d$
has
$\leq \!\ell d$
points
over
$\bbF_{\!\ell}$~\cite[the lemma in Chapter~1, Paragraph 5.2]{BS}.
For~a prime
ideal~$\frakl$
such that
$\ell \geq d+2$,
this is~contradictory.
}
\eop
\end{theo}

\section{Zariski density in the Hilbert scheme}

This~section is devoted to Zariski density of the counterexamples to the Hasse principle in the Hilbert~scheme. Our~result is, in fact, an application of the Zariski density in the moduli scheme established~above.

\begin{theo}
Let\/~$K$
be any number field,
$U_\reg \subset \Gr(2,15)_K$
the open subset of the Gra\ss mann scheme that parametrizes degree four del Pezzo surfaces, and\/
$\calHC_K \subset U_\reg(K)$
be the set of all degree four del Pezzo surfaces
over\/~$K$
that are counterexamples to the Hasse~principle.\smallskip

\noindent
Then\/~$\calHC_K$
is Zariski dense in\/
$\Gr(2,15)_K$.\medskip

\noindent
{\bf Proof.}
{\em
Let~us fix an algebraic closure
$\overline{K}$
and an embedding
of~$K$
into~$\overline{K}$.
Assume that, contrary to the assertion, 
$\calHC_K \subset U_\reg \subset \Gr(2,15)_K$
would not be Zariski~dense. It~is well-known that the Gra\ss mann scheme
$\Gr(2,15)_K$
on the right hand side is irreducible and projective of dimension
$(15-2) \!\cdot\! 2 = 26$.
The~subset
$\calHC_K$
must therefore be contained in a closed subscheme
$H \subset \Gr(2,15)_K$
of dimension at
most~$25$.

By Theorem~\ref{theo1}, the invariant map
$H \to \Pb(1,2,3)$
is dominant. Its generic fiber may thus be, possibly reducible, of dimension at
most~$23$.
In~particular, outside of a finite union of curves
$C \subset \Pb(1,2,3)$,
the special fibers are of dimension
$\leq \!23$,
as~well.

Now,~let us choose a
\mbox{$K$-rational}
point
$s \in [\Pb(1,2,3) \!\setminus\! C](K)$
that is the image of a degree four del Pezzo
surface~$S \in\calHC_K$
under the invariant~map. The~geometric fiber
$I^{-1}(s)_{\overline{K}}$
over~$s$
of the full invariant map
$\smash{I\colon U_\reg \to \Pb(1,2,3)}$
parametrizes all reembeddings of
$S$
into~$\smash{\Pb^4_{\overline{K}}}$
and is therefore a torsor under
$\PGL_5(\overline{K})/\Aut(S_{\overline{K}})$.
In particular, it is of
dimension~$24$.

This implies that 
$I^{-1}(s) \not\subseteq H$.
But~the orbit
of~$s$
under
$\PGL_5(K)$
para\-metrizes counterexamples to the Hasse principle, and is therefore contained
in~$H$.
As~$\PGL_5(K)$
is Zariski dense in
$\PGL_5(\overline{K})$,
this is a~contradiction.
}
\eop
\end{theo}

\setlength\parindent{0mm}
\end{document}